# The vehicle relocation problem

# for the one-way electric vehicle sharing

Maurizio Bruglieri[a], Alberto Colorni[a], Alessandro Luè[a,b]

[a]*Dipartimento INDACO - Politecnico di Milano, Via Durando 38/A, 20158 Milano, Italy*

[b]*Poliedra* - Politecnico di Milano, *via Giuseppe Colombo 40, 20133 Milano, Italy*

**Abstract**

Traditional car-sharing services are based on the two-way scheme, where the user picks up and returns the vehicle at the same parking station. Some services permits also one-way trips, which allows the user to return the vehicle in another station. The one-way scheme is quite more attractive for the users, but may pose a problem for the distribution of the vehicles, due to a possible unbalancing between the user demand and the availability of vehicles or free slots at the stations. Such a problem is more complicated in the case of electrical car sharing, where the travel range depends on the level of charge of the vehicles. The paper presents a new approach for the Electric Vehicle Relocation Problem, where cars are moved by personnel of the service operator to keep the system balanced. Such a problem generates a challenging pickup and delivery problem with new features that to the best of our knowledge never have been considered in the literature**.** We yield a Mixed Integer Linear Programming formulation and some valid inequalities to speed up its solution through a state-of-the art solver (CPLEX). We test our approach on verisimilar instances built on the Milan road network.

## 1. Introduction

The general expectation concerning personal mobility in western countries, confirmed by the investments of the principal car makers, is a shift from the internal combustion engine vehicles to electric ones. Such change will impact in particular on urban areas, which will experience less local emissions and a better air quality. Moreover, the electricity needed for the electrical vehicles will be generated more and more by clean and renewable technology such as hydroelectric, wind, geothermic and solar power plants. Another emerging element of the future mobility is idea of sharing the car among multiple users, so passing from the ownership of a vehicle to the use of a service.

At the light of such expectation, the *Green Move* project (Luè et al., 2012) has been studying and testing a new system based on small, electric and shared vehicles. Financed by Regione Lombardia, the project aims to face both the technology aspect and the service design part in order to identify a successful business model of vehicle-sharing. The basic idea is to create a flexible sharing service, based on electric vehicles and open to a wide range of different typology of users. The system will be made easily accessible thanks to a device, the Green-e-Box (Savaresi & Alli, 2011), which allows the inclusion into the vehicle sharing service of any vehicle. An institution is responsible for coordination, management, standardization and maintenance of a distributed system, while different organizations and companies own the subsystems (vehicles and docking stations). The interfaces of vehicles and docking stations are "standardized" and then fully inter-operable (i.e. each vehicle can make use of all the docking stations).

Each actor will be free to buy his own vehicles, on which he will impose his own restrictions and tariffs, or will contribute to the network by installing his own docking stations.

Within the *Green Move* project, the authors studied some optimization problems, related to the management of an electric car sharing service.

Car-sharing, understood as an organized form of shared use of the car, began to grow in Zurich in 1948 (Harms and Truffler, 1998). In the following years, particularly in the 70s and 80s, several car-sharing systems were set up, but still on a small scale, with poor results (Shaheen et al., 1998). The original idea has been gradually replaced by an offer structurally organized according to strict business criteria, in order to achieve economies of scale, which resulted in increased benefits to users in terms of low rates and diversification of the available fleet (Burlando and Manstretta, 2007). One of the main innovation consists in overcoming the two-way (or round trip) system, where the user picks up and return the vehicle at the same parking station. Some new services (such as Car2Go - www.car2go.com) permits also one-way trips, which allow the user to pick up the vehicle in one station, and return it in another one.

The design and management of a car-sharing service raise several optimization problems, which have been tackled in the literature (e.g. Du and Hall, 1997; George and Xia, 2011; Barth and Todd, 1999), in particular to determine the optimal size of the fleet and identify the location of the parking stations. For instance, Cucu et al. (2009) propose a methodology based on fuzzy logic algorithm, where the users' needs are modeled by way of suitable performance indicators, with the objective of ensuring the balance between costs, number of stations and level of service. The problems of the optimal location of the recharging stations and of the optimal electric vehicle routing are considered in Touati-Moungla and Jost (2010), where their mathematical programming formulations consider variants of well-known combinatorial optimization problems.

Traditional car-sharing services are based on the two-way scheme, where the user picks up and returns the vehicle at the same parking station. Some services permits also one-way trips, which allows the user to return the vehicle in another station

The one-way system is quite more attractive for the users, but may pose a problem for the relocation of the vehicles, due to a possible unbalancing between the demand and availability of vehicles (for example, near the railway stations at the beginning of a working day) or vice versa between the request for return of the vehicles and the availability of free slots. In such cases, the service provider has to develop strategies to reallocate the vehicles and restore an optimal distribution of the fleet of car-sharing service. Such a distribution could be based on the immediate needs at a particular station, or on an historical prediction (i.e. estimating the vehicle demand in the future in order to determine when and from where a relocation event occurs) (Barth and Todd, 1999). The activities of vehicle relocation can be carried out by the user itself or by the service provider (Barth et al., 2004). In the first case, the user is incentivized to car pool or to choose another location or reservation time; in the second case, which is the most common in the real services, the vehicles are physically transported using trucks or personnel. In the literature, also the platooning of the vehicles has been considered (Daviet et al., 1996), where the platoon is composed by a chain of technologically innovative vehicles, led by vehicle head. This technique, however, is still of little use, for safety reasons, and due to pending future technological developments (Chan, 2012).

In the following, a brief literature review regarding the methods proposed to solve the vehicle relocation problem is presented.

Chauvet et al. (1997) propose the use of a fleet of trucks, minimizing the number of cars to move between the stations and minimizing the travel time of the trucks. The first problem is an Uncapacitated Transportation Problem, solved through the Hitchcock algorithm (Hitchcock, 1941); while the second problem is addressed through a heuristic algorithm. One development is proposed in Chauvet et al. (1999), which makes the assumption of using a single auto transport truck that runs continuously on a circuit, arbitrarily established. Duron et al. (2000) presents a heuristics based on the immediate needs at the stations, i.e. the next station to be visited by the auto transport truck is chosen according to the current state of the system. In such a way, the algorithm gives priority to visit the stations that have the greatest likelihood of running out of vehicles. Di Febbraro (2012) represents the complex dynamics of the system, using a discrete event system simulation. The paper consider both the relocation made by both users and staff, simulating different scenarios, with the objective of reducing the number of required staff, and minimizing the number of car sharing vehicles to satisfy the system demand.

The problem of the relocation of *electric* vehicles (EVs) has been faced in Dror et al. (1998), which proposes an algorithm to manage auto transport trucks, based on a Tabu search approach and a savings heuristic (Toth and Vigo, 2002). The algorithm is applied to a car sharing service with fifty EVs and five stations, offered in the French town of Saint Quentin en Yvelines. The car sharing service offered in the same location was studied also by Hafez et al. (2001), which determined the needed number of auto transport trucks with an exact algorithm, and then minimized the total travel time of reallocation, studying three different heuristics.

**In our opinion the relocation approach based on auto transport trucks may be not well suitable for an urban** settings, from a practical point of view**, because stations may not be easily reachable by the trucks, and the operations of loading/unloading EVs is time consuming. For the EV relocation problem, we propose therefore** the **use of** a staff of car sharing operators (hereafter called *workers*). They may move easily and in eco-sustainable way from a delivery point to a pickup point using a folding bicycle that **can be loaded in the trunk of the EV which** need**s to be moved.** Such relocation approach generates a challenging pickup and delivery problem with features that, to the best of our knowledge, have been never considered in the literature. We call such a problem the Electric Vehicle Relocation Problem (EVRP). EVRP shares some features with the 1-*skip vehicle routing problem* (Archetti & Speranza, 2004) and the *rollon-rolloff problem* (Aringhieri et al., 2004; Bodin et al., 2000; De Muelemeester et al., 1997), i.e. the fact that just one item at the time can be picked up and delivered and routes starting and ending at a single depot cannot exceed a given maximum duration. However EVRP is more challenging than the above mentioned problems since it is complicated by the fact that the distance covered by a vehicle depends also by the item picked up, i.e. the residual electrical charge of the EV picked up. This further complication does not allow for instance to map the problem into a static bipartite graph like for the *rollon-rolloff problem* presented in Aringhieri et al. (2004) because the feasibility of an arc connecting a pickup request node with a delivery request node depends on the time when the pickup request node is reached since the residual charge of a parked EV increases over the time. In this work we yield the first Mixed Integer Linear Programming (MILP) formulation of the EVRP and we illustrate some methods to strengthen it and to speed up its solution through a state-of-the art solver (CPLEX). Finally we test our approach on verisimilar instances built on the Milan road network. The paper is organized as follows. In Section 2 we describe the problem, in Section 3 we present its MILP formulation, in Section 4 we introduce some methods to speed up its solution, in Section 5 we show some numerical results and in Section 6 we draw some conclusions.

## 2. Problem description

Let $L$ be the maximum distance that an EV can cover when its battery is fully charged. Such distance depends on the kind of EV considered; for instance $L$ can vary from 50 km for a Liberty Piaggio EV to 400 km for a Tesla EV (in the experimental campaign we assume that $L$ =150 km). Notice that when the battery of an EV is not fully charged, the maximum distance that can be covered is linearly proportional to the residual charge of the battery (i.e. an EV with residual charge at 50% can cover $L/2$ km). Concerning the recharge time $\Gamma$ of a battery, the question is slightly different since typically the recharge process comprises two phases: the first one is intensity-constant, the second one is tension-constant. The first phase allows to recharge the battery almost fully and it is linear on the time. The second phase is not linear on the time and can require some hours to achieve the full charge of the battery and to ensure an uniform recharge of all the cells that compose the battery. For sake of simplicity we do not consider the second phase of recharging to model the EVRP. The maximum time needed to complete the first phase depends on the recharge technology used and can vary for instance from $\Gamma = 1$ hour for a 380V Superfast Recharger to $\Gamma = 5$ hours for a 220V Multifast Recharger (in the experimental campaign we consider $\Gamma = 4$ hours).

Let $P$ be the set of pickup requests and $D$ the set of delivery requests. Each request $r \in P \cup D$ is characterized by a parking location $v_r$, i.e. a node of the road network, by the residual charge of the battery $\rho_r$, by a time $\tau_r$ which represents the earliest pickup time if $r \in P$, the latest delivery time if $r \in D$. Note that for a delivery request $r$, $\rho_r$ indicates the minimum charge level that the EV battery must have at time $\tau_r$ when the EV will be used, therefore if an EV is delivered before $\tau_r$ the charge level of its battery may be less than $\rho_r$ on condition that at least the charge level $\rho_r$ is achieved at the time $\tau_r$ recharging the battery after the delivery. We assume that the fleet of EV is homogeneous and therefore each delivery request can be satisfied picking up every EV of a pickup request on condition that it is compatible for time windows and battery charge level.

Given a team of $K$ workers which leave also at different times a single depot using folding bicycles, we want to determine their routes in such a way the duration of each route does not exceed a given threshold $T$ (i.e. the service time of the workers), each route ends in the depot, the number of requests served is maximized respecting the time windows and battery charge level constraints of each request.

## 3. Mathematical programming formulation

The formulation of the EVRP is inspired by the MILP formulation of the rollon-rolloff problem presented in Aringhieri et al (2004). We consider a graph $G=(N, A)$ that models all the possible actions rather than the road network. The set of nodes of $G$ is given by $N = P \cup D \cup \{0\}$ where 0 indicates the depot node. The set of arcs can be partitioned into two sets: the *EV arcs* and the *bike arcs*. The *EV arcs* model the action of a worker when he is moving by an EV from a pickup point to a delivery point; the *bike arcs* model the action of a worker when he is moving by bike from a delivery point or from the depot to a pickup point or to the depot. Therefore EV arcs $(i, j)$ are defined for each $i \in P$ and for each $j \in D$ such that $\tau_j \geq \tau_i + \frac{d_{ij}}{s'} + q' + q''$ and $d_{ij} \leq L$ where $d_{ij}$ indicates the minimum path to go from node $i$ to node $j$, $s'$ denotes the average speed of an EV, $q'$ is the time to park the EV and take the bike from the EV trunk and $q''$ is the time to load the bike in the EV trunk and exit with the EV from the parking. In similar way the bike arcs $(j,i)$ are defined for each $j \in D \cup \{0\}$ and for each $i \in P$ such that $\tau_i \geq \tau_j + \frac{d_{ij}}{s''} + q''$, where $s''$ denotes the average speed by bicycle.

The operational times $c$ associated with every kind of arcs are reported in Table 1.

| Arcs | Operational times | Involved nodes |
|------|-------------------|----------------|
| $(i, j)$ | $\frac{d_{ij}}{s'} + q' + q''$ | $\forall i \in P, \forall j \in D$ |
| $(j, i)$ | $\frac{d_{ji}}{s''}$ | $\forall i \in P, \forall j \in D$ |
| $(0,i)$ | $\frac{d_{0i}}{s''}$ | $\forall i \in P$ |
| $(j,0)$ | $\frac{d_{j0}}{s''}$ | $\forall j \in D$ |

**Table 1: Operational times $c$ of all arcs**

There are two main advantages to deal with the graph $G$ rather than directly with the road network. The first one is that to every feasible route of a worker corresponds always an elementary cycle on graph $G$, whereas this is not true in the original road network when there are multiple requests in the same parking and modeling non elementary cycles is by far harder (see Dror et al., 1998). The second advantage, even in the case of a single request for each parking, is that a formulation based on graph $G$ requires by far less variables than a formulation based on the road network, because variables are defined on the arcs and nodes of the used graph. The dimension of graph $G$ depends only by the number of requests ($|N|=|P|+|D|+1$ and $|A|<|N|^2$) and not by the number of the physical nodes (road intersections) and road links (e.g. the Milan road network considered in section 4 contains more than 23000 road links which are by far greater than $|A|$ even if 100 EVs need to be redistributed).

Let us introduce the binary routing variable $x_{ijk}$ equal to 1 if the $k$-th worker visits node $j \in N$ immediately after node $i$, 0 otherwise. Let us also introduce the continuous variables $t_{ik}$ to model the visit time of node $i$ by part of the $k$-th worker. We state that the EVRP can be modeled by way of the following MILP.

$$\max \quad \sum_{k=1}^{K} \sum_{(i,j) \in A: i \neq 0} x_{ijk} \tag{1}$$

subject to:

$$\sum_{j \in \delta^+(0)} x_{0jk} \leq 1 \quad \forall k = 1, ..., K \tag{2}$$

$$\sum_{k=1}^{K} \sum_{j \in \delta^+(i)} x_{ijk} \leq 1 \quad \forall i \in P \cup D \tag{3}$$

$$\sum_{j \in \delta^+(i)} x_{ijk} - \sum_{j \in \delta^-(i)} x_{jik} = 0 \quad \forall i \in P \cup D \cup \{0\}, \forall k = 1, ..., K \tag{4}$$

$$t_{ik} + c_{ij} x_{ijk} \leq t_{jk} + T(1 - x_{ijk}) \quad \forall (i,j) \in A: j \neq 0, \forall k = 1, ..., K \tag{5}$$

$$t_{ik} + c_{i0} x_{i0k} - t_{0k} \leq T \quad \forall i \in \delta^-(0), \forall k = 1, ..., K \tag{6}$$

$$t_{ik} \geq \tau_i \quad \forall i \in P, \forall k = 1, ..., K \tag{7}$$

$$t_{jk} \leq \tau_j \quad \forall j \in D, \forall k = 1, ..., K \tag{8}$$

$$d_{ij} x_{ijk} \leq L(\rho_i + \frac{t_{ik} - \tau_i}{\Gamma}) \quad \forall (i,j) \in A: i \in P, j \in D, \forall k = 1, ..., K \tag{9}$$

$$\rho_i + \frac{t_{ik} - \tau_i}{\Gamma} - \frac{d_{ij}}{L} x_{ijk} \geq \rho_j - \frac{\tau_j - t_{jk}}{\Gamma} - (\rho_j + 1)(1 - x_{ijk})$$
$$\forall (i,j) \in A: i \in P, j \in D, \forall k = 1, ..., K \tag{10}$$

$$1-\frac{d_{ij}}{L}x_{ijk} \geq \rho_j - \frac{\tau_j - t_{jk}}{\Gamma} - (\rho_j+1)(1-x_{ijk}) \quad \forall (i,j) \in A : i \in P, j \in D, \forall k = 1,...,K \quad (11)$$

$$x_{ijk} \in \{0,1\} \quad \forall (i,j) \in A, \forall k = 1,...,K \tag{12}$$

$$t_{ik} \geq 0 \quad \forall i \in P \cup D \cup \{0\}, \forall k = 1,...,K \tag{13}$$

Since each arc connects a pair of requests, the objective function represents the total number of requests served. Constraint (2) takes into account that at most $K$ workers are available and therefore at most $K$ routes can be generated imposing that at most one arc for each worker can leave the depot node 0. Constraints (3) impose that each request is satisfied at most once. Flow conservation constraints (4) ensure that the solution is a collection of cycles. Constraints (5) rule the time variables ensuring that the visit time of a node is given by the sum of the visit time of its predecessor and the operational time to go from the predecessor to the current node. Note that such constraints are not imposed for the depot node to ensure that the route can pass through the depot and at the same time they prevent the solution from containing isolated cycles that do not pass through the depot. In this way the formulation does not require additional subtour elimination constraints. Constraints (6) ensure that the duration of each route do not exceed the threshold $T$. Constraints (7) and (8) are the time windows respectively for the pickup requests and the delivery requests. Constraints (9) model the fact that the distance traveled by an EV is linearly proportional to the residual charge. Note that if $\rho_i + \frac{t_{ik} - \tau_i}{\Gamma} > 1$ such constraints become redundant since already the graph topology prevent the existence of arcs $(i, j)$ with $d_{ij} > L$. Finally constraints (10) and (11) ensure that an EV is delivered with a battery charge level such that at the time $\tau_j$ the charge level $\rho_j$ will be achieved.

### 4. Solving the MILP formulation of EVRP

In this section we present some methods to speed up the solution of MILP formulation (1)-(13) of the EVRP when a commercial MILP solver (e.g. CPLEX) is used.

First of all we note that the feasible region of MILP (1)-(13) may contain several optimal solutions since in an optimal solution, if any, the route of every worker can be exchanged with that one of any other worker yielding a different optimal solution in term of variables $x_{ijk}$ and $t_{ik}$. The presence of such multiple optimal solutions is harmful for a MILP solver since it may require more CPU time to detect the optimality of the solutions found. To prevent such situation we add to formulation (1)-(13) the following group of constraints that "breaks" the symmetry of the feasible region:

$$\sum_{(i,j) \in A : i \neq 0} c_{ij} x_{ijk'} \geq \sum_{(i,j) \in A : i \neq 0} c_{ij} x_{ijk''} \quad \forall k', k'' = 1,...,K : k' < k'' \tag{14}$$

Constraints (14) prevent the presence of multiple optimal solutions above mentioned since they impose that the routes are assigned to the workers according to the non- increasing operative cost ordering.

Another method to strengthen the MILP formulation of EVRP and speed up its solution consists in exploiting an upper bound $U$ on the maximum number of requests that can be satisfied. Indeed given $U$ we can strengthen the formulation adding to it the following valid inequality:

$$\sum_{k=1}^{K} \sum_{(i,j)\in A: i \neq 0} x_{ijk} \leq U \qquad (15)$$

Since the MILP formulation of EVRP can be solved in reasonable time for $K=1$, when $K > 1$ we can easily compute the upper bound $U$ solving with one worker the formulation (1)-(13) where constraints (7) and (8) are relaxed and the working time is extended to $KT$. The optimal value obtained in this way is an upper bound for the original problem because any feasible solution of the latter consists in at most $K$ routes, each one assigned to a worker, that satisfy all constraints (1)-(13) and therefore can be also traveled by one worker within the working time $KT$ satisfying certainly all constraints except (7) and (8). Hence any feasible solution of the original problem is also feasible for the formulation (1)-(13) where constraints (7) and (8) are relaxed and the working time is extended to $KT$.

A last trick to speed up the solution of the EVRP formulation with a MILP solver consists in providing it a starting feasible solution. To this purpose we consider the following simple but effective heuristic procedure based on the iterative solution with the MILP (1)-(13) of at most $K$ instances of the EVRP with one worker. At the $k$-th iteration the tour of the $k$-th worker is determined considering only the requests that are not already satisfied in the previous iterations until $k$ is equal to $K$ or all the requests are satisfied.

## 5. Experimental results

The section presents some numerical experiments made on the road network of Milano. The network is based on the database constructed by the Milan transport agency (Agenzia Mobilità Ambiente e Territorio, 2008), which contains information about the road geography, the nodes, permitted manoeuvres, and link attributes. The road network is made up by about 23000 road links and 32000 nodes. The experiments considers 9 parking stations located nearby some main attractors, i.e. Loreto, Cadorna, Porta Genova, Porta Garibaldi, Piola, Duomo, Stazione Centrale, Turati, San Babila (see Fig. 1).

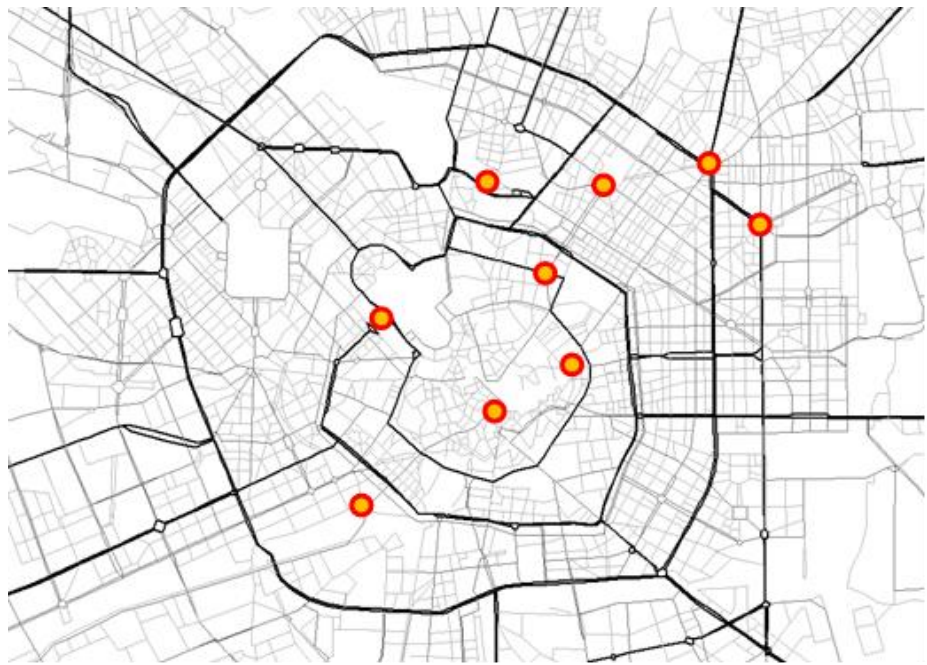

Fig. 1. The road network and the locations of the 9 charging stations considered in the numerical experiments.

We have built 60 instances of the EVRP considering randomly generated pickup and delivery requests on the 9 parking stations with random battery charge levels and random time windows between 8.00 a.m. and 3.00 p.m. In all instances $|P|=|D|$ and the total number of requests can assume the values 10, 20, 30, 40. For each value of $|P \cup D|$ we have built 5 instances of the EVRP and we have considered 3 values of $K = 1,2,3$.

The main input data values used for the MILP formulation presented in Sections 3 and 4 are summarized in Table 2.

| Input data | $T$ | $s'$ | $s''$ | $q'$ | $q''$ | $L$ | $\Gamma$ |
|---|---|---|---|---|---|---|---|
| Values | 300 minutes | 25 km/h | 15 km/h | 1 minute | 1 minute | 150 km | 240 minutes |

**Table 2: Main input data values used in the experiments**

The MILP formulation (1)-(13) has been implemented in AMPL (Fourer & Gay, 2002) and solved with the state of the art solver CPLEX11.0 on a PC Intel Xeon 2.80 GHz with 2GB RAM. The numerical results are reported in Table 3. The first column indicate the instance name, the second one the number of requests, the columns labeled "Served" report the percentage of requests satisfied, the columns labeled "CPU1" report the CPU time in seconds of MILP formulation (1)-(13) while the columns labeled "CPU2" report the CPU time in seconds obtained applying the speed up methods described in Section 4 (including also the CPU time spent to compute the upper bound $U$ and the starting feasible solution). The average results on each group of instances with the same request size are reported in Table 4. Here the columns labeled "Improv" report the average percentage relative improvement passing from CPU1 to CPU2.

Concerning to the quality of the solutions obtained it is possible to see that already one worker is able in average to satisfy the 70% of relocation requests and such percentage increases to 80% if we limit to instances up to 20 requests. With two workers the average percentage of relocation requests satisfied becomes the 80%. Whereas with three workers there is no improvement of requests satisfied for the instances up to 20 requests and the average percentage of relocation requests satisfied increases only of the 2%. Therefore 2 seems to be the suitable number of workers to satisfy instances up to 40 relocation requests on the test bed.

Concerning the CPU time the MILP formulation (1)-(13) is able to solve instances up to 20 requests in less than 1 seconds, but can require even more than 2 hours for instances with 40 requests. For the instances of such a size, the speed up methods presented in Section 4 reveals to be really useful. Indeed like shown in Table 4 they allow an average improvement of the CPU time of around 86% on the instances with 30 and 40 requests when the optimal number of workers $K$=2 is chosen.

| Instance | |P U D| | k=1 | | k=2 | | | k=3 | | |
|---|---|---|---|---|---|---|---|---|---|
| | | Served | CPU1 | Served | CPU1 | CPU2 | Served | CPU1 | CPU2 |
| Amat10_1 | 10 | 100% | 0.03 | 100% | 0.03 | 0.03 | 100% | 0.03 | 0.03 |
| Amat10_2 | 10 | 80% | 0.01 | 80% | 0.03 | 0.02 | 80% | 0.03 | 0.02 |
| Amat10_3 | 10 | 40% | 0.01 | 60% | 0.01 | 0.01 | 60% | 0.01 | 0.01 |
| Amat10_4 | 10 | 100% | 0.01 | 100% | 0.01 | 0.04 | 100% | 0.01 | 0.04 |
| Amat10_5 | 10 | 80% | 0.03 | 100% | 0.03 | 0.02 | 100% | 0.03 | 0.02 |
| Amat20_1 | 20 | 100% | 0.11 | 100% | 0.11 | 0.12 | 100% | 0.11 | 0.12 |
| Amat20_2 | 20 | 100% | 0.99 | 100% | 0.99 | 2.03 | 100% | 0.99 | 2.03 |
| Amat20_3 | 20 | 80% | 0.65 | 80% | 0.17 | 0.12 | 80% | 0.17 | 0.12 |
| Amat20_4 | 20 | 60% | 0.21 | 60% | 0.13 | 0.03 | 60% | 0.13 | 0.03 |
| Amat20_5 | 20 | 60% | 0.04 | 60% | 0.05 | 0.09 | 60% | 0.04 | 0.09 |
| Amat30_1 | 30 | 87% | 91.98 | 93% | 162.90 | 2.66 | 93% | 162.90 | 2.66 |
| Amat30_2 | 30 | 80% | 2.08 | 80% | 2.82 | 0.54 | 80% | 2.82 | 0.54 |
| Amat30_3 | 30 | 67% | 352.59 | 80% | 11.23 | 8.76 | 80% | 7.33 | 8.76 |
| Amat30_4 | 30 | 60% | 190.79 | 80% | 15.94 | 9.55 | 80% | 3.24 | 2.55 |
| Amat30_5 | 30 | 60% | 200.61 | 80% | 2506.17 | 358.12 | 87% | 1023.03 | 462.83 |
| Amat40_1 | 40 | 45% | 241.03 | 70% | 4030.65 | 651.96 | 75% | 8100.37 | 524.12 |
| Amat40_2 | 40 | 50% | 135.23 | 80% | 5440.26 | 263.11 | 80% | 2611.93 | 178.29 |
| Amat40_3 | 40 | 55% | 295.09 | 70% | 5388.41 | 831.85 | 75% | 6865.55 | 1544.58 |
| Amat40_4 | 40 | 50% | 467.67 | 65% | 4691.86 | 792.37 | 85% | 7221.15 | 581.22 |
| Amat40_5 | 40 | 50% | 343.04 | 65% | 3725.47 | 593.21 | 70% | 3984.21 | 872.64 |

**Table 3: Numerical results of EVRP formulation on AMAT data**

|  |P U D\| | k=1 | | k=2 | | | | k=3 | | | |
|---|---|---|---|---|---|---|---|---|---|---|
|  | Served | CPU1 | Served | CPU1 | CPU2 | Improv | Served | CPU1 | CPU2 | Improv |
| 10 | 80% | 0.02 | 88% | 296.41 | 0.02 | 99.99% | 88% | 0.02 | 0.02 | -9.09% |
| 20 | 80% | 0.40 | 80% | 0.29 | 0.48 | -64.83% | 80% | 0.29 | 0.48 | -65.97% |
| 30 | 71% | 167.61 | 83% | 539.81 | 75.93 | 85.93% | 84% | 239.86 | 95.47 | 60.20% |
| 40 | 50% | 296.41 | 70% | 4655.33 | 626.50 | 86.54% | 77% | 5756.64 | 740.17 | 87.14% |
| Average | 70% | 116.11 | 80% | 1372.961 | 175.731 | 52% | 82% | 1499.204 | 209.035 | 18% |

Table 4: Average numerical results of EVRP formulation on AMAT data

### 6. Conclusions

In this work a new approach to redistribute the vehicles of an EV sharing service has been proposed. It consists in moving the EVs by way of a team of workers that directly drive the EVs that need to be moved. The workers can easily move from a delivery point to a pickup point by way of a folding bicycle that can be loaded in the trunk of the EV that need to be moved. Such a problem generate a new challenging Vehicle Routing Problem for which we propose the first MILP formulation. Such a formulation is based on a graph representation of the problem rather than directly on the road network to the aim of avoiding non elementary cycles and to use less variables as possible. We test the formulation on big real world instances based on the road network of Milano (about 23000 road links and 32000 nodes). While for instances up to 20 requests the formulation solved by CPLEX 11.0 on a PC Intel Xeon 2.80 GHz with 2GB RAM requires less than 1 seconds, it can require even more than 2 hours for instances with 40 requests. To overcome such computational difficulties due to the NP-hardness of the EVRP we have contrived three expedients: we have added to the original formulation a group of constraints that "breaks" the symmetry of the feasible region; we have strengthened the formulation by way of a valid inequality based on an upper bound of the objective function; we have provided the MILP solver with a feasible starting solution obtained with a heuristic method based on the EVRP formulation.

Future work on the EVRP concerns the generation of the pickup and delivery requests in more verisimilar way exploiting the origin-destination traffic matrix yielded by the Milan transport agency. Moreover we have also an interest in investigating its combination with pricing policies i.e. the possibility of promoting parking stations with lack of EVs by decreasing the price of using an EV if it is delivered in such stations.

**Acknowledgements**

We thank Giovanni Alli and Andrea Giovanni Bianchessi of the Electronic and Information Technology Department of Politecnico di Milano for the technical information on EVs.